\documentclass[letterpaper, 10 pt, conference]{ieeeconf}  

\IEEEoverridecommandlockouts                              

\usepackage{amsmath} 
\usepackage{amssymb}  
\usepackage{mathtools}
\usepackage{algorithm}
\usepackage{algpseudocode}
\usepackage[para]{threeparttable}
\usepackage{booktabs}
\usepackage{mathtools}

\usepackage{array}
\newcolumntype{P}[1]{>{\centering\arraybackslash}p{#1}}
\newcolumntype{R}[1]{>{\raggedright\arraybackslash}p{#1}}

\newcommand{\bi}{\begin{itemize}} 
\newcommand{\ei}{\end{itemize}}
\newcommand{\be}{\begin{equation}} 
\newcommand{\ee}{\end{equation}}

\makeatletter
\let\NAT@parse\undefined
\makeatother
\usepackage[colorlinks=true,linkcolor=blue,citecolor=blue]{hyperref}
\usepackage[noadjust]{cite}

\title{\LARGE \bf
The Pickup and Delivery Problem with Crossdock for Perishable Goods
}

\author{Konstantinos Gkiotsalitis$^{1}$ and Amalia Nikolopoulou$^{2}$
\thanks{*This work was partially funded by the EU project CONDUCTOR}
\thanks{$^{1}$Dr Konstantinos Gkiotsalitis is with the School of Civil Engineering, Department of Transportation Planning and Engineering, National Technical University of Athens, Greece
        {\tt\small KGkiotsalitis@civil.ntua.gr}}%
\thanks{$^{2}$Dr Amalia Nikolopoulou is with the National Technical University of Athens, Greece
        {\tt\small ANikolopoulou@mail.ntua.gr}}%
}

\begin{document}

\maketitle
\thispagestyle{empty}
\pagestyle{empty}

\begin{abstract}
Our work departs from the original definition of the Pickup and Delivery Problem (PDP) and extends it by considering an interchange point (crossdock) where vehicles can exchange their goods with other vehicles in order to shorten their delivery routes and reduce their running times. Multiple operational constraints, such as time windows, vehicle capacities, and the synchronization of vehicles at the crossdock, are considered. In addition, the specific requirements of perishable goods, which should not be carried on long trips, are taken into account. Given this consideration, this study introduces the Pickup and Delivery Problem with Crossdock for Perishable Goods (PDPCDPG) and models it as a nonlinear programming problem. PDPCDPG is then reformulated to a MILP with the use of linearizations and its search space is tightened with the addition of valid inequalities that are employed when solving the problem to global optimality with Branch-and-Cut. Various computational experiments are conducted on benchmark instances found in the literature to assess the performance of our model. The results demonstrate the solution stability of the proposed approach.  The proposed model aims to provide a practical and effective approach for transportation and logistics companies dealing with time-sensitive deliveries. \newline


\end{abstract}

\section{INTRODUCTION}

The Pickup and Delivery Problem (PDP) is a well-known combinatorial optimization problem that plays a crucial role in transportation, logistics, and supply chain management \cite{ermis2021optimization,alesiani2022constrained}. Its objective is to determine efficient routes for vehicles to transport goods between pickup and delivery locations while satisfying various constraints such as vehicle capacity, time windows, and precedence relationships between tasks \cite{savelsbergh1995general}.

Traditional PDP models assume that all goods are non-perishable and can be transported without time constraints \cite{cordeau2010branch, xu2003solving}. However, in many real-world scenarios, there is a growing need to transport perishable goods that have specific time constraints and limited shelf life \cite{utama2020vehicle}. This gives rise to the Pickup and Delivery Problem with Perishable Goods (PDPPG), which introduces additional complexities and challenges.

The PDPPG requires careful consideration of the in-vehicle routing decisions to minimize the impact on the quality and freshness of the perishable goods. In particular, the in-vehicle ride times of perishable goods must be controlled to ensure they do not exceed their allowable shelf life. This poses a significant challenge for route planning and scheduling, as it involves optimizing the delivery performance while adhering to the perishability constraints. Similar issues arise in passenger delivery problems, such as the Dial-a-Ride Problem \cite{gkiotsalitis2021dial}, and the Problem of Synchronized Passenger Transfers at transfer stations and mobility hubs \cite{gkiotsalitis2018towards,gkiotsalitis2019robust,liu2021review,geurs2023smarthubs}.

In this paper, we aim to address the PDPPG by proposing a novel mathematical model and solution approach. Our approach considers the perishability constraints and optimizes the delivery performance by allowing vehicles to exchange goods at an interchange/transshipment point, also known as a crossdock. Unlike the classic PDP formulation where each vehicle starts and ends at the same depot, our formulation allows goods to have at most one transfer at the crossdock's location.

This extension of the problem, known as the Pickup and Delivery Problem with a Crossdock (PDPCD), offers opportunities to reduce vehicle running times and improve overall efficiency. By allowing goods to be transferred between vehicles, we can minimize the distance traveled and improve the utilization of available vehicles. Moreover, we explicitly consider the route duration constraints for pickup and delivery vehicles and the in-vehicle ride time constraints for perishable goods. This results in the introduction of the Pickup and Delivery Problem with Crossdock for Perishable Goods (PDPCDPG), which is formulated and analyzed in this study. Our contributions in this paper can be summarized as follows:

\begin{enumerate}
\item Formulation of the PDPCD with Perishable Goods which accounts for the in-vehicle ride time limitations of perishable goods with limited shelf life.
\item Modeling as a quadratic integer program and linearization.
\item Generation of test instances and performance evaluation assessing its stability and effectiveness in solving larger problem instances.
\end{enumerate}

Through these contributions, we aim to enhance the understanding and practical applicability of the PDPCDPG by providing a formal formulation, an optimized solution approach, and empirical analysis of its performance on diverse problem instances.

Our proposed PDPCDPG model aims to provide a practical and effective solution for transportation and logistics companies dealing with time-sensitive deliveries. By optimizing vehicle routing and considering the perishability constraints, we can achieve improved delivery performance and reduced costs.

The remainder of this paper is organized as follows. Section \ref{sec:2} provides a detailed review of related literature. Section \ref{sec:3} presents the mathematical formulation of our PDPCD model with perishable goods, including the specific constraints and objectives considered. It also includes the linearization of the model and related valid inequalities. Section \ref{sec:4} presents experimental results and performance analysis on a set of benchmark instances. Finally, section \ref{sec:5} concludes the paper and discusses potential future research directions.

\section{LITERATURE REVIEW}\label{sec:2}

Comparable to the Pickup and Delivery Problem with Crossdock is the Pickup and Delivery Problem with Transfers (PDPT), which deals with the transportation of goods. The initial definition is based on the general Pickup and Delivery Problem (PDP) \cite{savelsbergh1995general}, where packages are directly transported from suppliers to their corresponding customers by the same carrier. Berbeglia et al. \cite{berbeglia2007static} categorized PDPs based on the supplier-customer ratio, distinguishing between scenarios involving many suppliers delivering to many customers (many-to-many), a single supplier to a single customer (one-to-one), and the possibilities of shipments between these scenarios (one-to-many and others). The Pickup and Delivery Problems with intermediate facilities introduce the concept of an intermediate facility called a cross-dock (or transshipment point) where a vehicle can drop off or pick up loads. Unlike PDP, where the same vehicle fleet handles both pickup and delivery operations, PDPT relaxes the constraint that packages must be delivered by the same carrier, enabling package transfers between carriers. Mitrovic-Minic et al. \cite{mitrovic2006pickup} were among the first to develop a PDPT model with a single transfer point. Subsequent research explored solution techniques for PDPT with a single transfer node by \cite{cortes2010pickup,masson2013adaptive}, with Rais et al. \cite{rais2014new} allowing transfers to occur at any node. Masson et al. \cite{masson2014dial} introduced the concept of package transfers using shuttle services between two transfer points, while \cite{ghilas2016pickup} modeled passenger transportation with time windows and synchronized routes using scheduled lines for transfers.

Another closely related research area is the Vehicle Routing Problem with Cross-docking (VRPCD) for the transportation of goods. Wen et al. \cite{wen2009vehicle} laid the foundation for this area with their work on a practical distribution problem for a Danish company. They considered time window constraints for pickup and delivery nodes in the distribution network, as well as for the cross-dock (CD) facility, to account for a fixed planning horizon. Hasani et al. \cite{hasani2012capacitated} studied the VRPCD with split deliveries and multiple products, imposing simultaneous arrival of inbound vehicles at the CD, and formulating it as a mixed-integer linear programming (MILP) problem to optimize vehicle routes and the number of vehicles utilized. Tarantilis et al. \cite{tarantilis2013adaptive} proposed a multi-restart Tabu Search algorithm for solving benchmark instances from \cite{wen2009vehicle} and explored scenarios involving different vehicles for pickup and delivery operations, as well as total routing costs for open and closed route network configurations. Morais et al. \cite{morais2014iterated} developed three Iterated Local Search heuristic algorithms for VRPCD, testing them on instances from \cite{wen2009vehicle} and larger randomly generated instances with up to 500 nodes. Other works that approaches this problem with the use of heuristics or metaheuristics are \cite{dondo2013sweep,sadri2014modeling,birim2016vehicle,
nikolopoulou2017moving,grangier2017matheuristic,abad2018bi,baniamerian2019modified}.

Another research direction focuses on hybrid cross-docking distribution networks. Petersen et al. \cite{petersen2011pickup} investigated VRPCD with optional returns at the CD, and \cite{santos2013pickup} studied a hybrid network structure with cross-docking, allowing direct shipping between suppliers and customers. Nikolopoulou et al. \cite{nikolopoulou2019adaptive} examined a hybrid distribution network and used a Tabu Search algorithm to compare the transportation costs incurred by two alternative distribution strategies: direct shipping and cross-docking, considering transfers between two sets of origin and destination points. Gunawan et al. \cite{gunawan2021matheuristic} recently proposed a matheuristic algorithm for VRPCD without time windows for supplier and customer points. Guastaroba et al. \cite{guastaroba2016intermediate} provide a comprehensive review of distribution networks with freight consolidation and merging operations. Additionally, several works have specifically focused on the Pickup and Delivery Problem with Crossdock (PDPCD), introducing pickup and delivery routes to VRPCD (see \cite{petersen2011pickup,santos2013pickup}).

Studies related to perishable goods are mostly focused on the vehicle scheduling problem. Since Tarantilis 
and Kiranoudis \cite{tarantilis2001meta} introduced the Vehicle Scheduling Problem for Perishable Goods (VRPfPG), there have been several approaches that propose heuristics or metaheuristics for its solution given its NP-Hard nature \cite{hsu2007vehicle,osvald2008vehicle,chen2009production}. However, such studies do not consider crossdocking. 

The subsequent sections of this research paper focus on the static case of the Pickup and Delivery Problem with Crossdock for Perishable Goods (PDPCDPG) with  ride time limitations for goods. The goal is to minimize the running costs of the vehicles \cite{gkiotsalitis2023public}. 

\section{FORMULATION}\label{sec:3}
PDPCDPG is formally defined as follows. We consider a directed graph $\mathcal{G}=(\mathcal{V},\mathcal{A})$. The vertex set $\mathcal{V}$ is divided into three subsets: $\mathcal{O}\cup \mathcal{P}\cup \mathcal{D}$.

The set $\mathcal{O}$ consists of four copies of the depot, denoted as $\langle o_1,o_2,o_3,o_4\rangle$. These represent different stages of the vehicle's trip: $o_1$ represents the starting point of the vehicle for picking up goods, $o_2$ represents the return of the vehicle to the crossdock, $o_3$ represents the departure of the vehicle from the crossdock to deliver goods, and $o_4$ represents the end of the vehicle's trip after delivering all goods. It is important to note that the locations of $o_1$, $o_2$, $o_3$, and $o_4$ are the same since they all correspond to the depot location.  That is, the depot location is also the crossdock location.

The set $\mathcal{P}$ represents the pickup vertices, numbered $\langle 1,...,n\rangle$, and the set $\mathcal{D}$ represents the delivery vertices, numbered $\langle n+1,...,2n\rangle$. If two delivery requests have the same pickup location but different delivery locations, a duplicate pickup vertex is created. The same duplication process applies if two requests have the same delivery location but different pickup locations. However, if multiple goods share the same origin-destination pair, their requests can be represented by a single pickup and delivery pair. Therefore, each pickup and delivery vertex is associated with exactly one origin-destination pair. It follows that we have $n$ requests, where each request is a couple $(i,n+i)$ with $i\in \mathcal{P}$ being the pickup point and $n+i\in \mathcal{D}$ the associated delivery point for this origin-destination pair. 

The feasible arc set $\mathcal{A}$ is defined as follows: $\mathcal{A}=\big\{\{(o_1,j):j\in \mathcal{P}\}\cup\{(i,j):i\in \mathcal{P},j\in \mathcal{P},i\neq j\}\cup\{(i,o_2):i\in \mathcal{P}\}\cup \{(o_3,j):j\in \mathcal{D}\}\cup\{(i,j):i\in \mathcal{D}, j\in \mathcal{D},i\neq j\}\cup\{(i,o_4):i\in \mathcal{D}\}\big\}$. It is worth noting that a vehicle cannot directly travel from a pickup vertex to a delivery vertex without passing through the crossdock, represented by $o_2$ and $o_3$. This is a fundamental distinction from the Pickup-and-Delivery Problems with Transfers (PDPT).

In summary, a vehicle starts its trip from $o_1$, serves pickup vertices from the set $\mathcal{P}$, returns to the crossdock $o_2$ to exchange goods, departs from the crossdock $o_3$ to deliver the newly assigned goods to their delivery vertices, and finally returns to the depot $o_4$. 

To each vertex $i\in \mathcal{V}$, there is an associated pickup or delivery demand $q_i$ with $q_{i}\geq 0~\forall i\in \mathcal{P}$, and $q_{i}=q_{i-n}~\forall i\in \mathcal{D}$. This demand represents the number of goods of the origin-destination pair $(i,n+i)$. Note that vehicles start from $o_1$ empty and return to $o_4$ empty. That is, $q_{o_1}=q_{o_4}=0$. There is also a minimum service duration for boarding/alighting every product (good). If $\beta\in\mathbb{R}_{+}$ is the fixed time requirement for handling a single good, then this duration is $\beta \sum_{i\in \mathcal{P}}q_i$ at the crossdock location. For the crossdock location, we also assume an additional fixed time for unloading and reloading, $a\in\mathbb{R}_{+}$.

Let $\mathcal{K}$ be the set of vehicles. The capacity of vehicle $k\in \mathcal{K}$ is denoted as $Q_k\in\mathbb{R}_{+}$ and the maximum allowed duration of route $k$ as $T_k\in \mathbb{R}_{+}$. The cost and travel time of traversing a feasible arc $(i,j)~i\in \mathcal{A}$ without performing intermediate stops is $c_{ij}\in\mathbb{R}_{+}$ and $t_{ij}\in\mathbb{R}_{+}$, respectively. Note that the triangular inequality holds because both the costs and the travel times are non-negative. Let $L\in\mathbb{R}_{+}$ be the maximum allowed ride time of any perishable good due to its limited shelf life and $[e_i,l_i]$ the time window within which we should serve vertex $i$.

Let also $u_i^k\in\mathbb{R}_{+}$ be the time at which vehicle $k$ starts serving vertex $i\in \mathcal{P}\cup \mathcal{D}$ and $r_i\in\mathbb{R}_{+}$ the ride time of good $i$ corresponding to request $(i,n+1)$. We also introduce binary variables $x_{ij}^k,\eta_{i}^k$ and $\theta_i^k$. $x_{ij}^k$ is equal to 1 if vehicle $k$ serves vertices $(i,j)\in \mathcal{A}$ sequentially, i.e., vertex $j$ is served directly after vertex $i$. $\eta_i^k=1$ if vehicle $k$ unloads request $i\in \mathcal{P}$ to the crossdock. $\theta_i^k=1$ if vehicle $k$ reloads request $i\in \mathcal{P}$ from the crossdock. Binary variables $\tilde{\eta}_k$ and $\tilde{\theta}_k$ indicate also whether vehicle $k$ unloads or reloads at the crossdock, respectively. 
The problem has the following continuous, nonnegative variables, the first two of which have been already described: 
\begin{itemize}
\item $u_i^k$ is the time at which vehicle $k$ starts servicing vertex $i\in \mathcal{V}$ 
\item $r_i$ is the ride time of good $i$ corresponding to request $(i,n+i)$, where $i\in \mathcal{P}$
\item $\tau_k$ is the time at which vehicle $k\in \mathcal{K}$ finishes unloading at the crossdock
\item $w_k$ is the time at which vehicle $k\in \mathcal{K}$ starts reloading at the crossdock
\item $z_i$ is the time at which request $i\in \mathcal{P}$ is unloaded at the crossdock
\end{itemize}
The compact, three-index formulation of the PDPCDPG model is cast below. Note that we use a very large positive number $M\rightarrow +\infty$ for modeling purposes.

\begin{align}
& \min \sum_{k\in \mathcal{K}}\sum_{(i,j)\in \mathcal{A}}c_{ij}^k x_{ij}^k\label{eq1}\\
&\text{subject to:}\nonumber\\
&\sum_{k\in \mathcal{K}}\sum_{j:(i,j)\in \mathcal{A}}x_{ij}^k=1 ~~\forall i\in \mathcal{P}\cup \mathcal{D}\label{eq2}\\
& \sum_{i\in \mathcal{P}}\sum_{j:(i,j)\in \mathcal{A}}q_ix_{ij}^k\leq Q_k ~~\forall k\in \mathcal{K}\label{eq3}\\
& \sum_{i\in \mathcal{D}}\sum_{j:(i,j)\in \mathcal{A}}q_ix_{ij}^k\leq Q_k ~~\forall k\in \mathcal{K}\label{eq4}\\
&\sum_{j:(o_1,j)\in \mathcal{A}}x_{o_1j}^k=\sum_{j:(o_3,j)\in \mathcal{A}}x_{o_3j}^k=1 ~~\forall k\in \mathcal{K}\label{eq5}\\
&\sum_{j:(j,o_2)\in \mathcal{A}}x_{jo_2}^k=\sum_{j:(j,o_4)\in \mathcal{A}}x_{jo_4}^k=1 ~~\forall k\in \mathcal{K}\label{eq6}\\
&\sum_{i:(i,h)\in \mathcal{A}}x_{ih}^k-\sum_{j:(h,j)\in \mathcal{A}}x_{hj}^k=0 ~~\forall h\in \mathcal{P}\cup \mathcal{D},k\in \mathcal{K}\label{eq7}\\
&u_j^k\geq u_i^k+t_{ij}-M(1-x_{ij}^k) ~~\forall (i,j)\in \mathcal{A},k\in \mathcal{K} \label{eq8}\\
&e_i\leq u_i^k\leq l_i ~~\forall i\in \mathcal{V},k\in \mathcal{K} \label{eq9}\\
&\eta_i^k-\theta_i^k=\sum_{j\in \mathcal{P}\cup \{o_2\}:j\neq i}x_{ij}^k-\sum_{j\in \mathcal{D}\cup \{o_4\}:j\neq i+n}x_{i+n,j}^k \nonumber\\
&~~~~~~~~~~~~~~~~~~~~~~~~~~~~~~~~~~~~~~~~~~~~~ \forall i\in \mathcal{P},k\in \mathcal{K} \label{eq10}\\
&\eta_i^k+\theta_i^k\leq 1 ~~\forall i\in \mathcal{P},k\in \mathcal{K} \label{eq11}\\
&\frac{1}{M}\sum_{i\in \mathcal{P}}\eta_i^k\leq \tilde{\eta_k}\leq \sum_{i\in \mathcal{P}}\eta_i^k ~~\forall k\in \mathcal{K} \label{eq12}\\
&\frac{1}{M}\sum_{i\in \mathcal{P}}\theta_i^k\leq \tilde{\theta_k}\leq \sum_{i\in \mathcal{P}}\theta_i^k ~~\forall k\in \mathcal{K} \label{eq13}\\
&\tau_k=u_{o_2}^k+a\tilde{\eta}_k +\beta\sum_{i\in \mathcal{P}}q_i\eta_i^k ~~\forall k\in \mathcal{K}\label{eq14}\\
&w_k\geq \tau_k ~~\forall k\in \mathcal{K}\label{eq15}\\
&u_{o_3}^k=w_k+a\tilde{\theta}_k +\beta\sum_{i\in \mathcal{P}}q_i\theta_i^k ~~ \forall k\in \mathcal{K}\label{eq16}\\
&w_k\geq z_i-M(1-\theta_i^k) ~~ \forall i\in \mathcal{P},k\in \mathcal{K} \label{eq17}\\
&z_i\geq \tau_k-M(1-\eta_i^k) ~~ \forall i\in \mathcal{P},k\in \mathcal{K}\label{eq18}\\
&u_{o_2}^k-u_{o_1}^k\leq T_k ~~ \forall k\in \mathcal{K} \label{eq19}\\
&u_{o_4}^k-u_{o_3}^k\leq T_k ~~ \forall k\in \mathcal{K} \label{eq20}\\
&r_i= \sum_{k\in \mathcal{K}}\sum_{j:(j,n+i)\in \mathcal{A}}x_{j,n+i}^k u_{n+i}^k-\sum_{k\in \mathcal{K}}\sum_{j:(j,i)\in \mathcal{A}}x_{j,i}^k u_{i}^k \nonumber\\
&~~~~~~~~~~~~~~~~~~~~~~~~~~~~~~~~~~~~~~~~~~~~~~~~~~~~~~ \forall i\in \mathcal{P} \label{eq21}\\
&r_i\leq L ~~ \forall i\in \mathcal{P}\label{eq22}
\end{align}

The aim of the objective function \eqref{eq1} is to minimize the overall costs associated with operating the vehicles. The constraints outlined in \eqref{eq2} ensure that each vertex associated with a request is visited precisely once. To prevent exceeding the vehicle's capacity during the pickup and delivery process, constraints \eqref{eq3} and \eqref{eq4} are employed. In order to maintain a consistent structure, the pickup routes for each vehicle should commence from $o_1$, while the delivery routes should start from $o_3$, as indicated by constraints \eqref{eq5}. Similarly, constraints \eqref{eq6} guarantee that all pickup routes return to $o_2$, and all delivery routes return to $o_4$. Constraints \eqref{eq7} ensure that when a vehicle arrives at a pickup or delivery vertex, it must depart from that vertex to maintain flow conservation. To ensure proper sequencing, constraints \eqref{eq8} stipulate that if vehicle $k$ serves the arc $(i, j) \in \mathcal{A}$, then the starting time for servicing vertex $j$ by vehicle $k$ must be greater than or equal to $u_i^k$ plus the travel time from vertex $i$ to $j$. Moreover, constraints \eqref{eq9} guarantee that each vertex $i \in \mathcal{V}$ is serviced within its allotted time window. Let us now explain constraints \eqref{eq10}. When

\[\sum_{j\in \mathcal{P}\cup \{o_2\}:j\neq i}x_{ij}^k=1\]

vehicle $k$ will pickup request $i$. Similarly, when

\[\sum_{j\in \mathcal{D}\cup \{o_4\}:j\neq i+n}x_{i+n,j}^k=1\]

vehicle $k$ will deliver request $i$. Constraints \eqref{eq10} result in the following four cases: (a) if request $i$ is picked up but not delivered by vehicle $k$, then $\eta_i^k-\theta_i^k=1$ and because $\eta_i^k+\theta_i^k\leq 1$ (see constraints \eqref{eq11}) we have that $\eta_i^k=1$ and $\theta_i^k=0$; (b) if request $i$ is not picked up, but it is delivered by vehicle $k$, then $\eta_i^k-\theta_i^k=-1$ and because $\eta_i^k+\theta_i^k\leq 1$ we have that $\eta_i^k=0$ and $\theta_i^k=1$; (c) if request $i$ is not picked up and not delivered by vehicle $k$, then $\eta_i^k-\theta_i^k=0$ and because $\eta_i^k+\theta_i^k\leq 1$ we have that $\eta_i^k=0$ and $\theta_i^k=0$; (d) if request $i$ is picked up and delivered by vehicle $k$, then $\eta_i^k-\theta_i^k=0$ and because $\eta_i^k+\theta_i^k\leq 1$ we have that $\eta_i^k=0$ and $\theta_i^k=0$. 

Constraints \eqref{eq12} and \eqref{eq13} play a crucial role in determining whether vehicle $k$ unloads, reloads, or performs both actions at the crossdock. By considering the presence of unloading ($\tilde{\eta}_k$) and the total unloading requests ($\sum_{i\in\mathcal{P}}q_i\eta_i^k$), constraints \eqref{eq14} determine the completion time of unloading at the crossdock for vehicle $k$. To ensure a proper sequencing, constraints \eqref{eq15} guarantee that vehicle $k$ starts reloading at the crossdock only after it has completed the unloading process. Constraints \eqref{eq16} define the time at which vehicle $k$ finishes its unloading or reloading operations at the crossdock and is ready to depart. Constraints \eqref{eq17} ensure that reloading of request $i$ can only commence after it has been unloaded at the crossdock. Moreover, constraints \eqref{eq18} enforce that the unloading of request $i$ by vehicle $k$ must be completed before vehicle $k$ can finish its unloading process. By combining constraints \eqref{eq17} and \eqref{eq18}, it ensures that if a request is both unloaded and reloaded by different vehicles at the crossdock, the unloading vehicle must finish unloading before the reloading vehicle can reload the request. Constraints \eqref{eq19} and \eqref{eq20} maintain the maximum allowable duration $T_k$ for pickup routes from the trip start $o_1$ to the crossdock $o_2$ and delivery routes from the crossdock $o_3$ to the trip end $o_4$. Additionally, constraints \eqref{eq21} and \eqref{eq22} calculate the ride time for each request $i$ and ensure that it remains below the maximum allowable ride time $L$ of the perishable goods.

\subsection{Linearizations}

In order to linearize the nonlinear constraints \eqref{eq21}, we introduce a new variable $\tilde{u}_i$, which represents the time at which vertex $i\in \mathcal{P}\cup \mathcal{D}$ starts to be serviced. Since each vertex is served by exactly one vehicle, we can utilize this information to linearize the constraints \eqref{eq21}. The approach is to replace constraints \eqref{eq21} with the equality constraint $r_i=\tilde{u}_{n+i}-\tilde{u}_i,~ \forall i\in \mathcal{P}$ for all $i \in \mathcal{P}$. To achieve this, we need to enforce $\tilde{u}_i$ to take on the value of $u_i^{k^*}$, where $k^*$ represents the vehicle that serves vertex $i$. Introducing continuous slack variables $\sigma_{ij}^k$, we can express $\tilde{u}_i$ as follows:

\begin{equation} \label{eq23}
\begin{aligned}
& \tilde{u}_i+\sigma_{j,i}^k=u_i^k &&~~\forall (j,i)\in \mathcal{A}: i\in \mathcal{P}\cup \mathcal{D}, \forall k\in \mathcal{K}\\
&\sigma_{ji}^k\leq M(1-x_{j,i}^k)&&~~\forall (j,i)\in \mathcal{A}:i\in \mathcal{P}\cup \mathcal{D},k\in \mathcal{K}\\
&\sigma_{ji}^k\geq -M(1-x_{j,i}^k)&&~~\forall (j,i) \in \mathcal{A}: i \in \mathcal{P}\cup \mathcal{D}, k \in \mathcal{K}
\end{aligned}
\end{equation}

Constraints \eqref{eq24} together with constraints \eqref{eq23} can now replace the nonlinear constraints $\eqref{eq21}$:

\begin{align}
&r_i=\tilde{u}_{n+i}-\tilde{u}_i&\forall i\in \mathcal{P}\label{eq24}
\end{align}

\begin{figure}[H]
    \centering
    \includegraphics[scale=0.5]{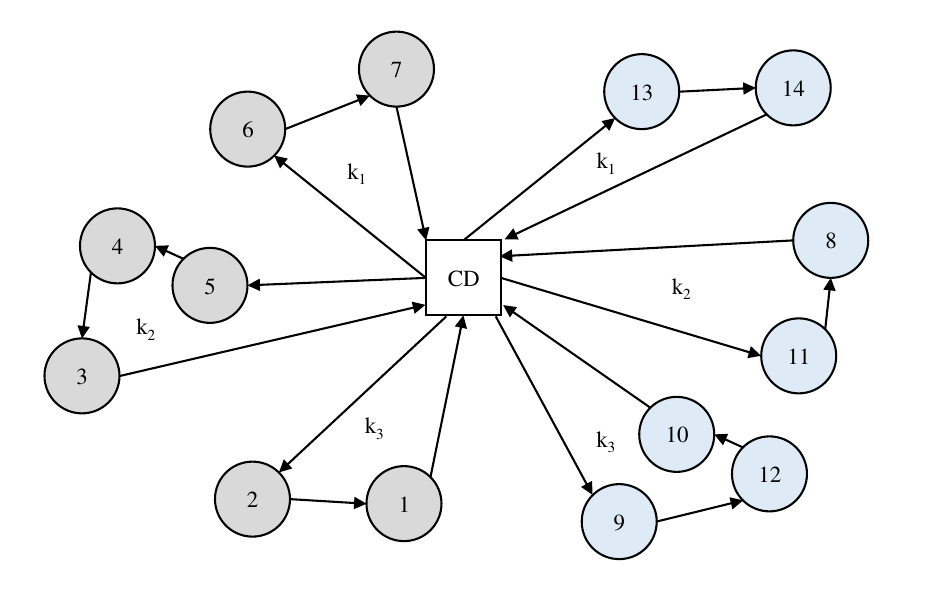}
    \caption{Illustrative example of a potential solution of the PDPCD with perishable goods in the case of three vehicles and seven requests.}
    \label{fig:1}
\end{figure}

To illustrate a potential solution of the PDPCDPG, Fig.\ref{fig:1} is provided. In this figure, 3 vehicles are assigned to 7 requests. Notice that vehicles 1-3 pick up the goods in locations 1-7, return to the crossdock where an exchange of goods takes place, and deliver them to their delivery points 8-14 before returning back to the depot, which is also the crossdock location.

\subsection{Valid Inequalities}
To enhance the {\em tightness} of our PDPCDPG formulation, we can introduce additional inequality constraints that do not eliminate any feasible solutions. These constraints should hold for any solution $\textbf{x}$ belonging to the feasible region $\mathcal{F}$ of our problem. By incorporating these valid inequalities, the computational time required to verify the feasibility of potential solutions increases due to the inclusion of extra inequality constraints. However, the tighter formulation allows us to exclude numerous non-optimal solutions, thereby facilitating a more focused search process. To tighten the formulation, we incorporate the following valid inequality constraints:

\subsubsection{Serve time tightening}
\begin{equation}
    u_i^k\geq e_i+\sum_{j:(j,i)\in A}\max\{0,e_j-e_i+t_{ij}\}x_{ji}~~\forall i\in P\cup D,k\in K
\end{equation}
\begin{equation}
    u_i^k\leq l_i+\sum_{j:(i,j)\in A}\max\{0,l_i-l_j+t_{ij}\}x_{ij}~~\forall i\in P\cup D,k\in K
\end{equation}
These valid inequalities were used in the past for solving the asymmetric Traveling Salesman Problem with time windows by branch and cut \cite{ascheuer2001solving}.

\subsubsection{Arc Elimination}
With arc elimination we remove arcs that are infeasible. In more detail:
\begin{itemize}
    \item arc $(i,j)\in A$ is infeasible if $e_i+t_{ij}>l_j$
    \item arcs $(i,j)$ and $(j',n+i)$ are both infeasible if $t_{i,j}+t_{j,o_2}+t_{o_3,j'}+t_{j',n+i}>L$ for $i\in P,j\in P,j'\in D$
\end{itemize}
\subsubsection{Sub-tour Elimination}
With sub-tour elimination we add valid inequalities that remove solutions which cannot be optimal, as follows:
\begin{itemize}
    \item $\sum_{k\in K}x_{ij}^k+\sum_{k\in K}x_{ji}^k\leq 1~~~~ \forall i\in P,j\in P$
    \item $\sum_{k\in K}x_{n+i,n+j}^k+\sum_{k\in K}x_{n+j,n+i}^k\leq 1~~~~ \forall i\in P,j\in P$
\end{itemize}
\subsubsection{Ride time: lower bound}
Finally, we add an additional set of valid inequalities related to the lower bounds of ride times to remove solutions that cannot be optimal:
\begin{equation}
r_i^k \geq t_{i,o_2}+(a\tilde{\theta}_k+\beta\sum_{i\in P}q_i\theta_i^k) +t_{o_3,n+i} ~~~~\forall i\in P,k\in K
\end{equation}

\section{Numerical Experiments}\label{sec:4}

\subsection{Demonstration in a toy network}

\begin{figure}
    \centering
    \includegraphics[scale=0.6]{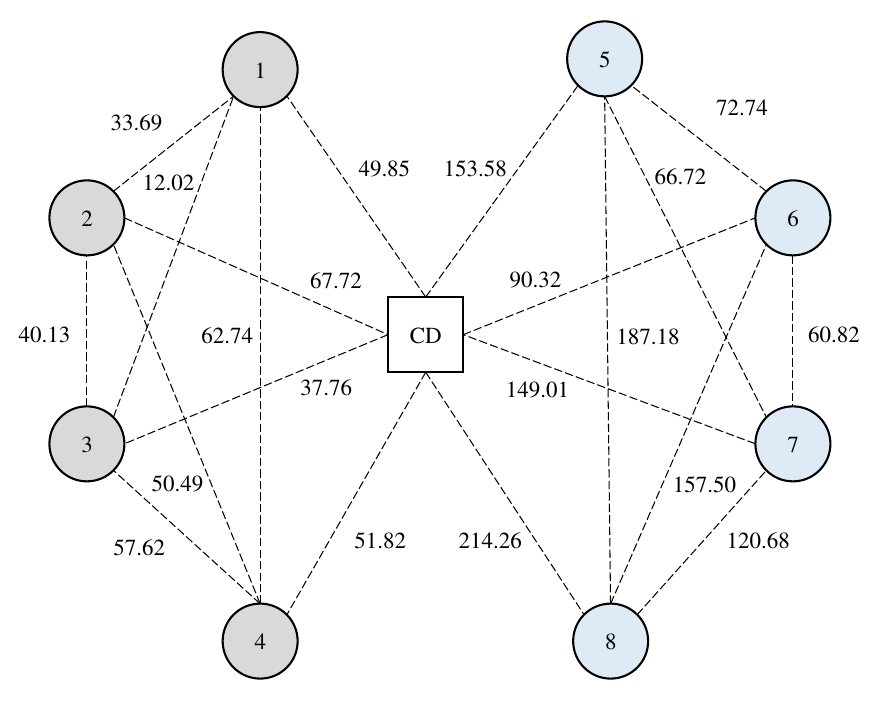}
    \caption{Network where vertices $o_1,o_2,o_3,o_4$ are at location 0 (depot). For visualization purposes, the presented arc travel times $t_{ij}$ are rounded to two decimal places.}
    \label{fig:2}
\end{figure}

We will now present the application of our PDPCD model for Perishable Goods in a toy network. The toy network, depicted in Fig.\ref{fig:2}, consists of 4 requests. The vehicle capacity is set to $Q_k=20$, and we have 2 vehicles available. The earliest possible start time for each vehicle is $e_{o_1}=360$, and the latest possible end time is $l_{o_4}=1320$. All time values are expressed in seconds.

At the crossdock, there is a fixed time of $a=10$ seconds allocated for unloading and reloading. Additionally, the handling time for a single product (good) at the crossdock is $\beta=1$ second.

For each vehicle ($k\in{1,2}$), the maximum allowed duration for a route is $T_k=480$ seconds. Furthermore, because of the limited shelf life of perishable goods, the ride time limit for any good is set to $L=550$ seconds.

The demand for pickup and delivery at each vertex is as follows: $q_1=q_5=16$, $q_2=q_6=10$, $q_3=q_7=4$, $q_4=q_8=4$.

In addition, Table \ref{tab:1} presents the time window for serving each vertex.

\begin{table}[H]
  \centering
  \small
  \caption{Lower and upper time for visiting any vertex $i\in P\cup D$}
    \begin{tabular}{rrrrr}
    \toprule
      $i$    & $e_i$  & $l_i$  & $e_{n+i}$ & $l_{n+i}$ \\
    \midrule
    1     & 442   & 562   & 823   & 943 \\
    2     & 455   & 575   & 852   & 972 \\
    3     & 360   & 471   & 793   & 913 \\
    4     & 475   & 595   & 1007  & 1127 \\
    \bottomrule
    \end{tabular}%
  \label{tab:1}%
\end{table}%

\begin{figure}
    \centering
    \includegraphics[scale=0.6]{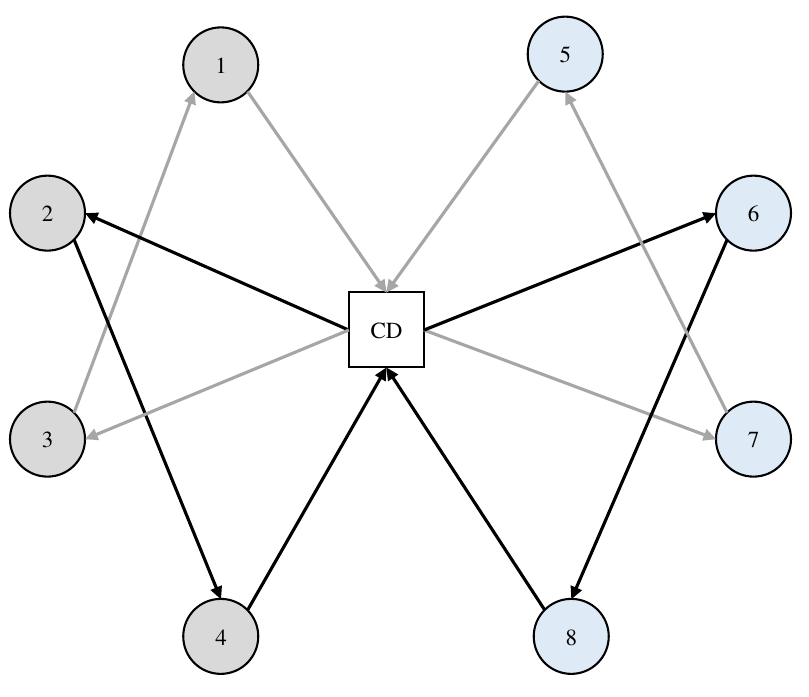}
    \caption{Optimal pickup and delivery routes of the two vehicles.}
    \label{fig:4}
\end{figure}

The optimal pickup and delivery routes of the two vehicles are presented in Fig.\ref{fig:4}, where gray color is used to illustrate the routes of the first and black color the routes of the second vehicle. These routes together with their associated travel times are summarized in Table \ref{tab:4}.

\begin{table}[htbp]
  \centering
  \caption{Travel time of each route and total travel cost for the service provider.}
    \begin{tabular}{llrr}
    \toprule
    vehicle & route type & served vertices & route travel time \\
    \midrule
    1     & pickup & $o_1 \rightarrow 3 \rightarrow 1  \rightarrow o_2$    & 99.813 \\
    2     & pickup & $o_1 \rightarrow 2 \rightarrow 4 \rightarrow o_2$    & 170.025 \\\\
    1     & delivery & $o_3 \rightarrow 7 \rightarrow 5  \rightarrow o_4$   & 369.310 \\
    2     & delivery & $o_3 \rightarrow 6 \rightarrow 8 \rightarrow o_4$   & 462.086 \\\midrule
      \multicolumn{3}{r}{total travel cost:}     & 1101.234 \\
    \bottomrule
    \end{tabular}%
  \label{tab:4}%
\end{table}%

Finally, the optimal start time of serving each vertex is presented in Table \ref{tab:5}. Note that these start times satisfy the imposed time windows. Table \ref{tab:5} reports also the ride time of each request. We note that:
\begin{itemize}
    \item the total travel cost of the vehicles is 1101.234 units.
    \item The ride times $r_i$ of the PDPCDPG solution are less than or equal to the maximum allowed ride time $L=550$, satisfying the ride time constraint for perishable goods. The same holds true for the total route travel times which are less than $T_k=480$.
\end{itemize}

\begin{table}
  \centering
  \small
  \caption{Start of service time $u_i$ at each vertex $i\in P\cup D$ and ride time $r_i$ for request $i\in P$.}
    \begin{tabular}{rrrrrrrrr}
    \toprule
      $i$    & $e_i$  & $u_i$ & $l_i$  & $e_{n+i}$ & $u_{n+i}$ & $l_{n+i}$ && $r_i$\\
    \midrule
    1     & 442  & {\bf 442.0} & 562   & 823  & {\bf 823.0} & 943 &&381.0\\
    2     & 455  & {\bf 544.5} & 575   & 852  & {\bf 852.0} & 972 &&307.5\\
    3     & 360  & {\bf 429.8} & 471   & 793  & {\bf 913.0} & 913 &&483.2\\
    4     & 475  & {\bf 595.0} & 595   & 1007 & {\bf 1009.5} & 1127 &&414.5\\
    \bottomrule
    \multicolumn{7}{l}{\footnotesize $u_{o_1}^1=u_{o_1}^2=360$}\\ 
    \multicolumn{7}{l}{\footnotesize $u_{o_2}^1= 599.99,u_{o_2}^2=646.82$}\\
    \multicolumn{7}{l}{\footnotesize $u_{o_3}^1= 643.99,u_{o_3}^2=743.77$} \\
    \multicolumn{7}{l}{\footnotesize $u_{o_4}^1= 1123.99,u_{o_4}^2=1223.77$ }
    \end{tabular}%
  \label{tab:5}%
\end{table}%

\subsection{Computational tests}

Herein, we present the results of our numerical experiments conducted on instances generated from the datasets introduced by Wen et al. \cite{wen2009vehicle} for the Vehicle Routing Problem with Crossdock (VRPCD). These datasets, originally released at \url{https://doi.org/10.11583/DTU.11786763.v1}, have been modified and expanded to accommodate the additional data requirements of our formulation, such as route duration, ride time constraints for perishable goods and a fixed number of vehicles. The vertex coordinates in the instances remain the same as those provided by Wen et al. \cite{wen2009vehicle}. 

Our branch and cut algorithm, which is employed to solve the MILP formulation of PDPCDPG, is implemented using Gurobi 9.0.3 in Python 3.7. The experiments were conducted on a server with a single thread and the following specifications: an Intel Xeon CPU E5-2650 v2 (2.60 GHz) processor and 16 GB of RAM. A time limit of four CPU hours was imposed for the execution of the algorithms. This allowed the exact branch and cut approach to solve instances with a maximum of 10 requests.

Due to the exponential growth in computational complexity for our NP-Hard problem, we restricted our numerical experiments to instances with a maximum of 10 requests. Table \ref{tab:8} presents the results obtained using the branch and cut method. The second column (CNS) indicates the number of constraints in each instance. The third column (NE) represents the number of explored nodes until the branch and cut algorithm's termination. The sixth column (ost) provides information on the performance of the globally optimal solution obtained by the branch and cut approach. Notably, the branch and cut method successfully solved instances with 10 requests within a time span of 3 CPU hours.

\begin{table}[H]
  \centering
  \begin{threeparttable}
  \caption{Results from Instances with up to 10 requests.}
  \footnotesize
    \begin{tabular}{lrrrrr}
    \toprule
          & &\multicolumn{3}{c}{BRANCH AND CUT} & \textbf{}   \\
    \cmidrule{3-5}
    {\bf Instance} & CNS & NE &  CPU (s) & ost  \\\midrule
    \textbf{4} & 1091  & 1        & 0.1   & 1101.23  \\
    \textbf{5} & 1536  & 39     & 0.4   & 1115.99 \\
    \textbf{6} & 3082  & 5,625  & 3.8   & 1390.48 \\
    \textbf{7} & 3988  & 9,892  & 13.7  & 1459.75 \\
    \textbf{8} & 5014  & 70,204  & 106.1 & 1556.51 \\
    \textbf{9} & 6160  & 213,507  & 1196.3 & 1518.77\\
    \textbf{10} & 9891  & 2,611,554  & 10917.7 & 1823.83 \\
    \bottomrule
    \end{tabular}%
    \begin{tablenotes}
      \small
      \item \scriptsize {\bf CNS}: problem constraints, {\bf NE}: nodes explored by Branch and Cut, {\bf ost}: cost of the globally optimal solution.
    \end{tablenotes}
  \label{tab:8}%
  \end{threeparttable}
\end{table}%

\section{Conclusion}\label{sec:5}

In this research, we proposed a nonlinear model for the multi-vehicle Pickup and Delivery Problem with Crossdocking for Perishable Goods. To solve this problem, we linearized the model and formulated it as a mixed-integer linear programming problem. Our experimental results demonstrate that the proposed PDPCDPG model can be solved to global optimality within a reasonable time for instances with up to 10 demand requests.

Our study extended the original pickup and delivery problem formulation by incorporating a crossdock, which can effectively reduce vehicle running costs by allowing goods to change vehicles, and ride time limitations due to the limited shelf life of perishable goods. In future research, further investigation can be conducted to consider the use of multiple interchange points, leading to a more generalized formulation with multiple transfers. Furthermore, additional operational constraints, such as the heterogeneity of the fleet of vehicles in terms of capacities and costs, could be considered.

Finally, given the NP-Hard nature of the PDPCDPG problem, which is inherited from the NP-Hardness of the PDPCD problem \cite{santos2013pickup}, future research can explore problem-specific heuristics or metaheuristics that can produce (sub)optimal solutions for larger problem instances with more than 10 requests.

\section*{ACKNOWLEDGMENT}

This work was partially funded by the European Union's Horizon Europe research and innovation programme CONDUCTOR (Grant Agreement no 101077049).

\bibliographystyle{IEEEtran}
\bibliography{IEEEabrv,citationlist}

\end{document}